\documentclass[12pt]{article}
\usepackage[utf8]{inputenc}
\usepackage[french,english]{babel}
\usepackage[T1]{fontenc}
\usepackage{amsmath}
\usepackage{amsfonts}
\usepackage{amssymb}
\usepackage{graphicx}
\usepackage[dvipsnames]{xcolor}
\usepackage{tcolorbox}
\usepackage{hyperref}
\usepackage[nottoc]{tocbibind}
\usepackage[right=2.5cm, left=2.5cm, top=2.5cm, bottom=2.5cm]{geometry}
\usepackage{cite}
\usepackage{setspace}
\newtheorem{Theorem}{Theorem}[section]
\newtheorem{Lemma}{Lemma}[section]
\newtheorem{Proposition}{Proposition}[section]
\newtheorem{Definition}{Definition}[section]
\newtheorem{Remark}{Remark}[section]
\newtheorem{Corollary}{Corollary}[section]

\newgeometry{includefoot,left=2cm,right=2cm,bottom=1cm,top=2cm}
\title{ \LARGE \bf%
{$L^p$-asymptotic stability of 1D damped wave equations with nonlinear and localized damping}
\thanks{This research was partially supported by the iCODE Institute, research project of the IDEX Paris-Saclay, and by the Hadamard Mathematics LabEx (LMH) through the grant number ANR-11-LABX-0056-LMH in the ``Programme des Investissements d'Avenir''.}
}
\author{Meryem Kafnemer$^{1}$, and Yacine Chitour$^{2}$}
\footnotetext[1]{UMA, ENSTA Paris, Institut Polytechnique de Paris, F-91120 Palaiseau, France}
\footnotetext[2]{L2S, Universit\'e Paris Saclay, France.}
\begin{document}
	\maketitle
	\begin{abstract}
		In this paper, we study the $L^p$-asymptotic stability with $p\in (1,\infty)$ of the one dimensional nonlinear damped
		wave equation with a localized damping and Dirichlet boundary conditions in a bounded domain $(0,1)$. We start by addressing the well-posedness problem. We prove the existence and the uniqueness of weak solutions for $p\in [2,\infty)$ and the existence and the uniqueness of strong solutions for all $p\in [1,\infty)$.  The proofs rely on the well-posedness already proved in the $L^\infty$ framework by \cite{Chitour2019} combined with a density argument. Then we prove that the energy of strong solutions decays exponentially to zero. The proof relies on the multiplier method combined with the work that has been done in the linear case in \cite{kafnemer2022}.
	\end{abstract}
\section{Introduction}
In this paper we consider the damped one dimension wave equation with nonlinear and localized damping and Dirichlet boundary conditions which yields  the following nonlinear problem:
\begin{equation}\label{probnl}
\left\lbrace
\begin{array}{cccc}
z_{tt} -z_{xx} + a(x)g(z_t)=0 & \hbox{for } (t,x) \in  \mathbb{R_+} \times (0,1), \\
z(t,0)= z(t,1)=0   & \hbox{for } t\geq 0 ,\\
z(0,\cdot)= z_0\ ,\ z_t(0,\cdot)=z_1 & \hbox{on }\  (0,1),
\end{array}
\right.
\end{equation}
where initial data $(z_0,z_1)$ belong to an $L^p$ functional space to be defined later, where $p\in [1,\infty)$. The function $a$ is a continuous non-negative 
function on $[0,1]$, bounded from below by a positive constant on some non-empty open interval $\omega$ 
of $(0,1)$, which represents the region of the domain where the damping term is active. The nonlinearity $g: \mathbb{R} \mapsto \mathbb{R}$ is a $C^1$ non-decreasing function such that $g(0)=0$, $g'(0)> 0$, and
$
g(x)x\geq 0 \ \ \forall x\in \mathbb{R}.$ 
\\ \\
The nonlinear problem \eqref{probnl} has already been studied several times in a hilbertian framework, i.e., with $p=2$. The well-posedness is a classical result and is a consequence of the theory of maximal monotone operators  (see for instance \cite{Martinez2000} for the non-localized case and \cite{Kafnemer2021} for the localized case). Exponential stability was treated using the multiplier method that was generally presented in \cite{Komornik1994} and then used in different contexts of the linear problem, we refer to \cite{Alabau} and \cite{Martinez1999} for more details and extensive references in the Hilbertian framework. The same method was adopted and used in the nonlinear problem by \cite{Martinez2000} to prove the exponential stability in two dimensions with no localization and with a growth hypothesis on the nonlinearity $g$. Then the exponential stability for the localized case was established in \cite{Kafnemer2021} with the same hypotheses.\\ \\
However, when it comes to the non-hilbertian framework, fewer results exist since it is still an unusual framework for the study of PDEs and has only been considered recently. As a result, we have less tools and techniques available to work with in such framework. The usual well-posedness proof techniques based on Hilbert Spaces (maximal monotone operators for instance) are no longer usable in this framework. Similarly, the multiplier method as known in the hilbertian framework to prove exponential stability does not work and requires at least to be generalized of the multipliers in order to achieve stability results.
\\ \\
The main results that exist and are relevant to our context are primarily gathered in \cite{Chitour2019}, \cite{Haraux2009} and \cite{kafnemer2022}.
 	The $p$-th 
 	energy $E_p$ of a solution that has been used in the three references has been introduced first in \cite{Haraux2009} as a generalization of the standard 
 	Hilbertian energy $E_2$. It is an equivalent energy to the natural energy in the $L^p$ framework and is defined by:
 	\begin{equation}\label{energyp}
 	E_p(t)=\frac{1}{p}\int_{\Omega}\left( |z_x(t,x)+z_t(t,x)|^p + |z_x(t,x)-z_t(t,x)|^p\right)dx.
 	\end{equation} 
 	The reference \cite{Haraux2009} also provides some useful energy estimates that were used for instance in \cite{Amadori2020},\cite{Chitour2019} and \cite{kafnemer2022} to obtain stability results.
 	\\ \\
 In the nonlinear case, \cite{Haraux2009} proves the well-posedness for all $p\geq 2$ using an argument based on the well-posedness in the Hilbertian framework and an equivalent energy functional but with global growth conditions on $g$.
 Useful $L^p$ estimates have also been provided in \cite{Haraux2009}, which led to proving polynomial decay of the energy in the nonlinear case, with $g$ a non-decreasing $C^1$ function behaving like $ks|s|^r$, $r,k>0$. We extend in this paper the results of \cite{Haraux2009} by proposing a well-posedness proof for weak solutions for all $2 \leq p< \infty$ and a proof for strong solutions for all $1 \leq p< \infty$. The proof is based on a density argument combined with the well-posedness already established in $L^\infty$ by \cite{Chitour2019}. We also extend the stability result by proving an exponential stability with no additional growth hypotheses on the nonlinearity $g$.
 \\ 
 \\ 
 Always in the nonlinear case but with a linearly bounded damping, an argument based on D'Alembert formula and fixed point theory is used in \cite{Chitour2019} to prove the well-posedness for all $p\geq2$. The reference proves the existence of solutions in $L^\infty$ framework for any nonlinear $g$ satisfying the hypotheses mentioned in the beginning of the introduction. However, for the well-posedness in an $L^p$ framework, the damping is supposed to be (uniformly) linearly bounded to be able to use the fixed point argument with D'Alembert formula. The latter reference relies on Lyapunov 
techniques for linear time varying systems along with estimates inspired from \cite{Haraux2009} to prove $L^p$ semi-global exponential stability in the nonlinear problem under restrictive hypotheses on initial data (imposed to belong to $L^\infty$ functional spaces) and for $ p\geq 2$; other stability results have been shown in the same reference in particular $L^\infty$ stability but always with more conditions on initial data. Still in \cite{Chitour2019}, They another stability result is obtained when the nonlinearity is (uniformly) linearly bounded using an interpolation for the initial data: semi-global exponential stability for $q$ satisfying $2 \leq q < p$ with initial data belonging to both $L^p$ spaces and $L^2$ is established. We extend the well-posedness results in $L^\infty$ of \cite{Chitour2019} to $L^p$ frameworks for $p\geq 2$ and also their well-posedness results in $L^p$ for $p\geq 2$  by removing the assumption of a linearly bounded nonlinearity. The latter well-posedness result is also stated in our work for $1\leq p<2$ for strong solutions. Additionally, we extend the stability results of \cite{Chitour2019} by providing a semi-global exponential stability result for strong solutions with no additional restrictions on initial data or the nonlinearity $g$.
\\ \\
The paper can also be seen as an extension to the nonlinear case for $1<p<\infty$ of \cite{Kafnemer2021} where the well-posedness was established for the linear problem for all $p>1$ using the argument based on D'Alembert formula from \cite{Chitour2019} and 
an exponential decay in the linear problem for all $1<p<+\infty$ was proved using a generalized multiplier method. The same reference also establishes an exponential stability when $p=\infty$ and $p=1$ in some cases of a constant global damping. 
\\ 
\\
We use the work that has already been done in the linear case in \cite{kafnemer2022} alongside with some techniques from \cite{Martinez2000} and \cite{Chitour2019} to provide a proof for the exponential stability of strong solutions in the nonlinear case. The proof is based on a linearizing principle
 to reduce the study of the nonlinear problem to that of the linear problem. 
\\ \\ 
The paper is organized as follows: in Section \ref{statement}  we properly state the problem with the functional framework and the considered hypotheses, we also rewrite the problem using Riemann invariants and we prove that the $p$-th 
 	energy $E_p$ is non-increasing. In Section \ref{sect4.3}, we give the proof of the well-posedness of the problem  for $1\leq p<\infty$ by proving the existence and the uniqueness of weak solutions for $2\leq p<\infty$ and the existence and the uniqueness of strong solutions for $1\leq p<\infty$. In Section \ref{sect4.4} we prove the exponential decay  of the energy by treating first an auxiliary linear problem and then by concluding the result for the nonlinear case.
\section{Statement of the problem}\label{statement}
	Consider Problem~\eqref{probnl} where we assume the following hypotheses satisfied:\\ \\
	$\mathbf{ (H_1)}$ $ a : [0,1] \rightarrow \mathbb{R}\ $ is a non-negative continuous function such that
	\begin{equation} \label{a0}
	\exists \ a_0 >0, \ a \geq a_0\ \ \hbox{on}\ \ \omega= (b,c) \subset [0,1],
	\end{equation}
	where $\omega$ is a non empty interval such that $b=0$ or $c=1$, i.e., $\bar{\omega}$ contains a
	neighborhood of $0$ or $1$. There is no loss of generality in assuming $c=1$, taking $0$ as an observation point.
	\\ \\
$\mathbf{ (H_2)}$
$g: \mathbb{R} \to \mathbb{R}$ a $C^1$ non-decreasing function such that $g(0)=0$, $g'(0)> 0$, and
\begin{align}\label{gxx}
g(x)x\geq 0 \ \ \forall x\in \mathbb{R},
\end{align}
We define now the functional framework. For $p \in [1, \infty]$, consider the function spaces
\begin{align}
X_p&:= W^{1,p}_0(0,1) \times L^p(0,1), \\
Y_p&:= \left(W^{2,p}(0,1) \cap W^{1,p}_0(0,1)\right) \times W^{1,p}_0(0,1),
\end{align}
where $X_p$ is equipped with the norm
\begin{align}
\Vert{(u,v)}\Vert_{X_p}:&= \left(\frac{1}{p}\int_0^1 \left(|u' + v|^p + |u' - v|^p\right) dx \right)^\frac{1}{p},   \ \ \text{ if } 1\leq p <\infty. \\
\Vert{(u,v)}\Vert_{X_\infty}:&= ||u' + v||_\infty + ||u' - v||_\infty, \ \ \text{ if }  p=\infty,
\end{align}
and the space $Y_p$ is equipped with the norm
\begin{align}
\Vert(u,v)\Vert_{Y_p}:&=
\left(\frac{1}{p}\int_0^1 \left(|u'' + v'|^p + |u'' - v'|^p\right) dx \right)^\frac{1}{p}, \ \ \text{ if } 1\leq p <\infty.\\
\Vert{(u,v)}\Vert_{X_\infty}:&=||u' + v||_\infty + ||u' - v||_\infty, \ \ \text{ if }  p=\infty.
\end{align}
Initial conditions $(z_0,z_1)$ for weak (resp.\ strong) solutions of \eqref{probnl} are taken in $ X_p$ (resp.
in $ Y_p$). 
\begin{Definition} The solutions of \eqref{probnl} are defined as follows.
	
	\begin{description}
		\item[$(i)$] For all $(z_0,z_1) \in X_p$, the function $$z\in L^\infty(\mathbb{R_+}, W^{1,p}_0(0,1)) \cap W^{1,\infty}(\mathbb{R_+},L^p(0,1)) 
		$$ is said to be a weak solution of Problem \eqref{probnl} if it satisfies the problem in the dual sense (meaning that the equalities are taken in the weak topology of $X_p$). 
		\item\item[$(ii)$] For all $(z_0,z_1) \in Y_p$, the function $$z\in L^\infty(\mathbb{R_+},W^{2,p}(0,1)\cap W^{1,p}_0(0,1))\cap W^{1,\infty}(\mathbb{R_+},W^{1,p}_0(0,1))$$ is said to be a strong solution of Problem \eqref{probnl} if it satisfies the problem in the classical sense.
	\end{description}
\end{Definition}
\ \\
We define the Riemann invariants for all $(t,x) \in \mathbb{R_+}\times (0,1)$ by
\begin{align}\label{Riemann}
\rho(t,x)= {z_x(t,x)+z_t(t,x)}, \\
\xi(t,x)= {z_x(t,x)-z_t(t,x)}.
\end{align}
Along strong solutions of \eqref{probnl}, we deduce that
\begin{equation}\label{probw}
\left\lbrace
\begin{array}{ll}
\rho_t - \rho_x  =-a(x)g\left(\frac{\rho-\xi}{2}\right) &\text{in }   \mathbb{R_+} \times (0,1), \\
\xi_t + \xi_x  =a(x)g\left(\frac{\rho-\xi}{2}\right) &\text{in }   \mathbb{R_+} \times (0,1), \\
\rho(t,0)-\xi(t,0)= \rho(t,1)-\xi(t,1)=0   &  \forall t \in \mathbb{R_+} ,\\
\rho_0:=\rho(0,.)= {y_0'+ y_1}\ ,\ \xi_0 :=\xi(0,.)={y_0'-y_1},
\end{array}
\right.
\end{equation}
with $\left(\rho_0, \xi_0 \right) \in W^{1,p}(0,1)\times W^{1,p}(0,1)$.
\\ \\
The $p$th-energy $E_p$ of a solution $z$  defined in \eqref{energyp} can also be written as
\begin{align}
E_p(t)= \frac{1}{p} \int_0^1 (|{\rho}|^p+|{\xi}|^p) dx.
\end{align}
We prove the same proposition as Proposition (2.2) from \cite{kafnemer2022}. 
\begin{Proposition}\label{prop1}
	Let $p \in [1, \infty)$ and suppose that a strong solution $y$ of \eqref{probnl} exists and is defined on
	a non trivial interval $I \subset \mathbb{R_+}$ containing $0$, for some initial conditions $(z_0,z_1)\in
	Y_p$. For $t\in I$, define
	\begin{align}\label{phidef}
	\Phi (t) := \int_0^1 [{\cal{F}}(\rho)+{\cal{F}}(\xi)] dx,
	\end{align}
	where $\rho$, $\xi$ are defined in \eqref{Riemann} and ${\cal{F}} $ is a $C^1$ convex function. Then $\Phi$ is well defined
	for $t\in I$ and satisfies
	\begin{align}\label{phidec}
	\frac{d}{dt} \Phi(t) =-\int_0^1 a(x) g\left(\frac{\rho-\xi}{2}\right)({\cal{F}}'(\rho)- {\cal{F}}'(\xi)) dx \leq 0.
	\end{align}
\end{Proposition}
\textbf{\emph{Proof.}} The proof is similar to the proof of \cite[Proposition 2.2]{kafnemer2022}. We obtain by following the same steps that
\begin{align}
\frac{d}{dt} \int_0^1 ({\cal{F}}(\rho)+{\cal{F}}(\xi)) dx
=- \int_0^1 a(x)g\left(\frac{\rho-\xi}{2}\right)({\cal{F}}'(\rho)- {\cal{F}}'(\xi)) dx. 
\end{align}
Thanks to the convexity of $\cal{F}$ and the hypothesis on $g$ that states that $g(x).x >0 $ for all $x\neq0$, we conclude that
\begin{align}
-\int_0^1 a(x)g\left(\frac{\rho-\xi}{2}\right)({\cal{F}}'(\rho)- {\cal{F}}'(\xi)) dx \leq 0,
\end{align}
which concludes the proposition.	\begin{small}
	\begin{flushright}
		$\blacksquare$
	\end{flushright}
\end{small}
\begin{Remark}
	The previous proposition has been first introduced in \cite{Haraux}  and reused in \cite{Chitour2019} to prove that the energy functional is non-increasing. Then it was improved in \cite{kafnemer2022} by omitting the hypothesis that function $\mathcal{F}$ should be even on top of being convex. 
\end{Remark}
	Before we state the next result, we introduce for $r\geq 0$ the following notation
	\begin{align}
	\lfloor x\rceil^r:=\textrm{sgn}(x)|x|^r ,\ \   \forall x \in \mathbb{R},
	\end{align}
	where $\textrm{sgn}(x)=\frac{x}{\vert x\vert}$ for nonzero $x\in\mathbb{R}$ and {$\textrm{sgn}(0)=[-1,1]$}.
	We have the following obvious formulas which will be repeatedly used later on:
	\begin{align}
	\frac{d}{dx}(\lfloor x\rceil^r)=r|x|^{r-1}, \ \ \forall r\geq 1,\ x\in\mathbb{R},\\
	\frac{d}{dx}(|x|^r)=r\lfloor x \rceil^{r-1}, \ \ \forall r> 1,\ x\in\mathbb{R}.
	\end{align}
	\begin{Corollary}\label{cor:decrease}
For $(z_0,z_1) \in Y_p$, one has that along strong solutions and for $t\geq 0$,
		\begin{align}\label{E_p'}
		E_p'(t) = -\int_0^1 a(x)g\left(\frac{\rho - \xi}{2}\right)\left(\lfloor \rho \rceil^{p-1}-
		\lfloor \xi \rceil^{p-1} \right) dx,
		\end{align}
		where $\rho$ and $\xi$ are the Riemann invariants defined in \eqref{Riemann}. Moreover, for $(z_0,z_1) \in X_p$, suppose that the solution $z$ of \eqref{probnl}
		exists on $\mathbb{R_+}$. Then the energy $t\longmapsto E_p(t)$ is non-increasing.  
	\end{Corollary}
	The first part of the corollary is an immediate application of Proposition \ref{prop1}
	while the second part is obtained by a standard density argument.

\section{Well-posedness} \label{sect4.3}
Before adressing stabilization issues, we start by studying the well-posedness of the problem. One should note that our proof of the well-posedness in the linear problem is no longer applicable to the nonlinear problem. The reason is that we would find ourselves with a fixed-point problem for a map that does not necessarily map to itself. Hence, the requirement of a new proof for the nonlinear case. \\
\begin{Theorem}\label{wellposedness}
Suppose Hypotheses $\mathbf{ (H_1)}$ and $\mathbf{ (H_2)}$ are satisfied, then for all initial conditions $(z_0,z_1)\in X_p$ with  $2\leq p <\infty $, we have the existence of a unique weak solution $z$ such that
	\begin{align}
	z\in L^\infty(\mathbb{R_+}, W^{1,p}_0(0,1)) \text{ and } z_t\in L^\infty(\mathbb{R_+},L^p(0,1)).
	\end{align}
	Moreover, if $(z_0,z_1)\in Y_p$ with 	 $1\leq p <\infty $ then we have the existence of a unique strong solution $z$ such that
	\begin{align}
	z\in L^\infty(\mathbb{R_+},W^{2,p}(0,1)\cap W^{1,p}_0(0,1)) \text{ and }  z_t\in L^\infty(\mathbb{R_+},W^{1,p}_0(0,1)).
	\end{align}
\end{Theorem}
\ \\
\emph{\textbf{Proof:}}\\ \\
\textbf{Weak solutions:}
Fix $2\leq p< +\infty$ and let $(z_0,z_1)\in X_p$. 
\\ \\ 
Since $X_\infty$ is dense in $X_p$ for all $2\leq p<\infty$ and $Z_0=(z_0,z_1)\in X_p$, there exists a sequence $\lbrace Z_0^n \rbrace_n \subset X_\infty$ such that $Z_0^n \to Z_0$ in $X_p$.
\\ \\ 
Since $Z_0^n \in X_\infty$, we have thanks to \cite[Theorem 1]{Chitour2019} the existence of a unique solution $Z^n=(z^n,z_t^n)$ such that $(z^n,z_t^n)\in L^\infty(\mathbb{R_+};W^{1,\infty}_0(0,1))\times W^{1,\infty}(\mathbb{R_+};L^\infty(0,1))$. Moreover we have that
for all $t\geq 0$
\begin{align}
||(z^n,z^n_t)||_{X_\infty} \leq 2 \max\left(||{z^n_0}'||_{L^\infty(0,1)}, ||z^n_1||_{L^\infty(0,1)}  \right).
\end{align}
We prove now that, for every $t_0\geq 0$, the sequence $\lbrace Z^n(t_0,\cdot) \rbrace_n$ is a Cauchy sequence in $X_p$. Define for all $(t,x)\in \mathbb{R_+}\times (0,1)$ the quantity $e^{n,m}$ as 
\begin{align}
e^{n,m}=z^n-z^m,
\end{align}
which is a solution of the following problem 
\begin{equation}\label{probe}
\left\lbrace
\begin{array}{cccc}
e^{n,m}_{tt} -e^{n,m}_{xx} + a(x)\left(g(z^n)-g(z^m)\right)=0 & \hbox{for } (t,x) \in  \mathbb{R_+} \times (0,1), \\
e^{n,m}(t,0)= e^{n,m}(t,1)=0   &  t\geq 0 ,\\
e^{n,m}(0,\cdot)= z_0^n-z_0^m\ ,\ e^{n,m}_t(0,\cdot)=z_1^n-z_1^m.
\end{array}
\right.
\end{equation}
The energy of $e^{n,m}$ at time $t$ is denoted by $E_p(e^{n,m})(t)$ and is non-increasing. Indeed, if we use the Riemann invariants of Problem \eqref{probe} denoted by $\rho(e^{n,m}),\ \xi(e^{n,m})$ given by $\rho(e^{n,m})= \frac{e^{n,m}_x+e^{n,m}_t}{2}$ and $\xi(e^{n,m})= \frac{e^{n,m}_x-e^{n,m}_t}{2}$ and rewrite the problem like Problem \eqref{probw}, then we follow the same steps in Proposition \ref{prop1} and Corollary \eqref{cor:decrease}, we obtain that along strong solutions of \eqref{probe}
\begin{align}
E_p(e^{n,m})'(t) =-2\int_0^1 \left(g(z^n_t)-g(z^m_t)\right)\left({\cal{F}}'(\rho(e^{n,m}))- {\cal{F}}'(\xi(e^{n,m})) \right) dx,
\end{align}
where $\mathcal{F}$ is taken to be the convex function $ | \cdot |^{p}$. By simple manipulations, it follows that 
\begin{align}
&E_p(e^{n,m})'(t) \notag 
\\&=-2\int_0^1 \frac{\left(g(z^n_t)-g(z^m_t)\right)}{z^n_t-z^m_t}\left({\cal{F}}'(\rho(e^{n,m}))- {\cal{F}}'(\xi(e^{n,m})) \right)(\rho(e^{n,m})-\xi(e^{n,m})) dx \leq 0,
\end{align}
which confirms using a density argument that $E_p(e^{n,m})$ is non-increasing along weak solutions. It follows that 
\begin{align}
E_p(z^n-z^m)(t)=E_p(e^{n,m})(t) \leq E_p(e^{n,m})(0)=E_p(z_0^n-z_0^m),
\end{align}
which gives that for every $t_0\geq 0$, $\lbrace Z^n(t_0,\cdot) \rbrace_n $ is a Cauchy sequence in $X_p$ since $\lbrace Z_0^n \rbrace_n$ is a Cauchy sequence in $X_p$. It follows then that $\lbrace Z^n(t_0,\cdot) \rbrace_n $ converges to a limit in $X_p$ that we denote by $Z(t_0,\cdot)=(z(t_0,\cdot),z_t(t_0,\cdot))$, in particular for $t_0=0$ we have that $Z(0,\cdot)=Z_0(\cdot)$. Note also that the convergence is uniform with respect to $t_0\geq 0$. We define the function $(t,x)\mapsto Z(t,x)$, where $Z(t_0,\cdot)$ is the limit of $Z^n(t_0,\cdot)$ in $X_p$ for every $t_0\geq 0$. 
\\ \\ 
We need to prove now that the limit $Z$ is a weak solution of \eqref{probnl}. Fix $T>0$ and denote $\psi$ a test function that belongs to $C^1([0,T]\times[0,1])$
 also verifying $\psi(T,\cdot)=\psi(0,\cdot)\equiv 0$ and {$\psi(\cdot,0)=\psi(\cdot,1)\equiv 0$}. Define
\begin{align}
A_T^n(\psi)=	\int_0^T \int_0^1 (z_{tt}^n-z_{xx}^n)\psi \, dxdt.
\end{align}
We have that
\begin{align}
A_T^n(\psi)-A_T^m(\psi)= \int_0^T \int_0^1 (z_{tt}^n-z_{xx}^n)\psi \, dxdt -\int_0^T \int_0^1 (z_{tt}^m-z_{xx}^m)\psi \, dxdt,
\end{align}
by integrating by part, it follows that
\[
A_T^n(\psi)-A_T^m(\psi)= -\int_0^T \int_0^1 (z_{t}^n-z_{t}^m)\psi_t \, dxdt +\int_0^T \int_0^1 (z_{x}^m-z_{x}^n)\psi_x \, dxdt.
\]
Using Holder's inequality,
\begin{align}
\left|A_T^n(\psi)-A_T^m(\psi) \right| &\leq  \left(\int_0^T \int_0^1|z_t^n-z_t^m|^p dxdt\right)^{\frac{1}{p}}   \left(\int_0^T \int_0^1|\psi_t|^q dxdt\right)^{\frac{1}{q}}
\notag 
\\&+\left(\int_0^T \int_0^1 |z_x^n-z_x^m|^p dxdt\right)^{\frac{1}{p}}   \left(\int_0^T \int_0^1 |\psi_x|^q dxdt\right)^{\frac{1}{q}},
\end{align}
which means that 
\begin{align}
\left| A_T^n(\psi)-A_T^m(\psi) \right| \leq T^{\frac{1}{p}} E_p(e^{n,m})^{\frac{1}{p}}(0) \left(||\psi_t||_{L^q((0,T)\times (0,1))}+||\psi_x||_{L^q((0,T)\times (0,1))}\right).
\end{align}
By a density argument we obtain that for all $\psi$ in the space
\begin{align}
	\mathcal{X}_q^T=&\lbrace \psi\ : \ [0,T]\times [0,1] \mapsto \mathbb{R_+} \ : \  \ (\psi,\psi_t)\in W^{1,q}((0,T)\times(0,1)) \times L^q((0,T)\times(0,1)) ,
	\notag\\
	&\ \ \ \ \ \ \ \ \ \ \ \ \ \ \ \psi(T,\cdot)=\psi(0,\cdot)\equiv 0 \text{ and } \psi(\cdot,0)=\psi(\cdot,1)\equiv 0\rbrace,
\end{align} 
where $q$ is the conjugate exponent of $p$ and is equal to $\frac{p}{p-1}$ for $p\geq 2$
we have that
\begin{align}
\left| A_T^n(\psi)-A_T^m(\psi) \right| \leq T^{\frac{1}{p}} E_p(e^{n,m})^{\frac{1}{p}}(0) \left(||\psi_t||_{L^q((0,T)\times (0,1))}+||\psi_x||_{L^q((0,T)\times (0,1))}\right),
\end{align}
and then 
\begin{align}
\left| A_T^n(\psi)-A_T^m(\psi) \right| \leq T^{\frac{1}{p}} E_p(e^{n,m})^{\frac{1}{p}}(0) ||\psi||_{\mathcal{X}_q^T},
\end{align}
which gives that $\lbrace A_T^n \rbrace_n $ is a Cauchy sequence in $(\mathcal{X}_q^T)'$, the dual of $\mathcal{X}_q^T$, since $\lbrace Z_0^n\rbrace_n$ is in $X_p$. We conclude then that $\lbrace A_T^n \rbrace_n $ converges in $(\mathcal{X}_q^T)'$, i.e.,  $\lbrace  z^n_{tt} -z^n_{xx}\rbrace_n$ converges to ${z}_{tt} -{{z}}_{xx}$ as a linear functional on $\mathcal{X}_q^T$.
\\ \\ 
We also know that for all $n\in \mathbb{N}$, 
\begin{align}
z^n_{tt} -z_{xx}^n= -g(z^n_t), \text{ on } \mathbb{R_+}\times (0,1),
\end{align}
which means that, the sequence of linear functionals $\lbrace -g(z^n_t) \rbrace_n$  defined on $\mathcal{X}_q^T$ also converges ${z}_{tt}(t,\cdot)-{{z}}_{xx}(t,\cdot)$ in $(\mathcal{X}_q^T)'$. 
\\ \\
We use now the existence of weak solutions in $L^2$ framework (see \cite{Martinez2000}, \cite{Kafnemer2021}). For  $(z_0,z_1)\in X_p$ with $p\geq 2$, we have the existence of a unique weak solution ${z}\in W^{1,\infty}(\mathbb{R_+},L^2(0,1))\cap L^\infty(\mathbb{R_+},H^1_0(0,1))$. The solution $z$ satisfies that for almost every $t\in \mathbb{R_+}$
\begin{align}
{z}_{tt}(t,\cdot)-{z}_{xx}(t,\cdot)=-g({z}_t(t,\cdot) \ \ \ \text{ in } H^{-1}(0,1).
\end{align}
In particular, for every $T>0$, ${z}_{tt}(t,\cdot)-{z}_{xx}(t,\cdot)=-g({z}_t)$ belongs to $(\mathcal{X}_2^T)'$. Since $q\leq 2$, one has that  $(\mathcal{X}_q^T)'\subset (\mathcal{X}_2^T)'$, yielding in particular that  $z$ is a weak solution of System \eqref{probnl} in $X_p$.
\\
\\ \\
\textbf{Strong solutions:}
Take $Z_0=(z_0,z_1)\in Y_p$ for all $1\leq p<\infty$, using Sobolev embeddings classical results, we have that $Z_0=(z_0,z_1)\in X_\infty$.
\\ \\ 
We use \cite[Theorem 1]{Chitour2019} to have the existence of a unique solution $Z^n=(z,z_t)$ such that $$(z,z_t)\in L^\infty(\mathbb{R_+}; W^{1,\infty}_0(0,1) )\times W^{1,\infty}(\mathbb{R_+};L^{\infty}(0,1)).$$ Moreover we have the following inequality also proved in \cite[Theorem 1]{Chitour2019} for all $t\geq 0$
\begin{align}
||(z,z_t)||_{X_\infty} \leq 2 \max\left(||z_0'||_{L^\infty(0,1)}, ||z_1||_{L^\infty(0,1)}  \right).
\end{align}
We are going to use Proposition \ref{prop1} for $w=z_t$, where $z$ is a solution of Problem \eqref{probnl}.  By differentiating \eqref{probnl} with respect to $t$, we obtain that $w$ satisfies the following problem
\begin{equation}\label{weq2}
\left\lbrace
\begin{array}{ll}
w_{tt} - w_{xx}  =-4w_t g'(w) &\text{in }   \mathbb{R_+} \times (0,1), \\
w(t,0)=w(t,1)=0   &  \forall t \in \mathbb{R_+} ,\\
w(0,\cdot)=  z_1,\  w_t(0,\cdot)=z_0''-g(z_1).
\end{array}
\right.
\end{equation}
We define for all $(t,x)\in \mathbb{R_+}\times (0,1)$ the Riemann invariants for \eqref{weq2},
\begin{align}
u=w_x+w_t,\\
v=w_x-w_t.
\end{align}
Along strong solutions of \eqref{weq2}, we have 
\begin{equation}\label{weqr}
\left\lbrace
\begin{array}{ll}
u_{t} - u_{x}  =-2(u-v) g'(w) &\text{in }   \mathbb{R_+} \times (0,1), \\
v_{t} + v_{x}  =2(u-v) g'(w) &\text{in }   \mathbb{R_+} \times (0,1), \\
u(t,0)-v(t,0)=u(t,1)-v(t,1)=0   &  \forall t \in \mathbb{R_+}.
\end{array}
\right.
\end{equation}
We define the $p$-th energy associated with $w$ as
\begin{align}
E_p(w)(t)=\frac{1}{p}\int_0^1(|u|^p+|v|^p)\, dx.
\end{align}
By Corollary \eqref{cor:decrease}, $E_p(w)$ is non-increasing along solutions $w=z_t$, which implies that 
\begin{align}
E_p(w)(t) \leq E_p(w) (0),\quad \textrm{ for a.e. } t\geq 0.
\end{align}
Then, using the fact that $E_p(w)$ is an equivalent energy to $\frac{1}{p}\int_0^1(|w_x|^p+|w_t|^p)\, dx$, it follows that 
\begin{align}
\int_0^1|w_x|^p\, dx \leq CpE_p(w) (0),
\end{align}
which means that
\begin{align}\label{w}
||z_t||_{W^{1,p}(0,1)} \leq  (CpE_p(w) (0))^{\frac{1}{p}}.
\end{align}
This implies that $(z,z_t)\in Y_p$ for all $t\geq 0$ which yields the required regularity for a strong solution.
\begin{small}
	\begin{flushright}
		$\blacksquare$
	\end{flushright}
\end{small}
 \begin{Remark}
 	For strong solutions in the case $p\geq 2$, we can easily use the results that have been proved in \cite{Haraux2009} for $p\geq 2$ to prove the well-posedness. Indeed, let $(z_0,z_1)\in Y_p$, with $ p\geq 2$ this implies that $(z_0,z_1)\in \left(H^2(0,1)\cap H^1_0(0,1)\right)\times H^1_0(0,1)$. We have then the existence of a unique strong solution $$z\in C(\mathbb{R_+},H^1_0(0,1))\cap C^1(\mathbb{R_+},L^2(0,1)),$$ such that $$(z(t,\cdot),z_t(t,\cdot) )\in \left(H^2(0,1)\cap H^1_0(0,1)\right)\times H^1_0(0,1), \ \ \forall t\in \mathbb{R_+},$$ which means that $z_t(t,\cdot)\in L^\infty(0,1)$ we can then use \cite[Corollary 2.3, item (ii) ]{Haraux2009} that implies in our context that if $(z_0,z_1)\in Y_p$ then the solution $z\in L^\infty(\mathbb{R_+},W^{2,p}(0,1)\cap W^{1,p}_0(0,1))$ and $z_t\in L^\infty(\mathbb{R_+},W^{1,p}_0(0,1))$ which guarantees the well-posedness in $Y_p$.
 \end{Remark}
\begin{Remark}
	The reason why the argument based on D'Alembert formula and fixed point theory that was used in \cite{Chitour2019} cannot be used in the nonlinear case without imposing the extra assumption that $g$ is linearly bounded, is that the fixed point argument cannot be used when we cannot prove that the nonlinearity maps the convex compact on which we define the fixed point formula to itself.
\end{Remark}
\section{Exponential stability}\label{sect4.4}
In this section we are interested in the asymptotic stability of Problem \eqref{probnl}. Our goal is to prove that the energy along strong solutions of Problem \eqref{probnl} is exponentially decreasing, with an exponential rate of decrease depending on the $Y_p$-norm of the initial data. (This property is usually referred as semi-global exponential stability in the control literature.) To do so, we plan to use the work that has already been done in the linear case to treat the nonlinear  case. The main stability result that we achieve is given by the following theorem.
\begin{Theorem}\label{stability}
	Assume $\mathbf{(H_1)}$ and $\mathbf{(H_2)}$ satisfied and $1<p<\infty$. Given $(z_0,z_1) \in Y_p$, there exists a constant $C_p(z_0,z_1)>0$ that depends on the norm of initial conditions in $Y_p$ such that for all $t\in \mathbb{R_+}$
	\begin{align}\label{expdec}
	E_p(t) \leq E_p(0)e^{1-C_p(z_0,z_1)t}.
	\end{align}
\end{Theorem}
\begin{Remark}
	The stability rate in Theorem \ref{stability}  depends on initial conditions. That is due to the fact that we are only manipulating strong solutions in the proof and we need an estimate of the $W^{1,p}(0,1)$ norm of $z_t(t,\cdot)$.
\end{Remark}
\subsection{Asymptotic stability of an auxiliary linear problem}
To prove the exponential stability of Problem \eqref{probnl}, we are going to start by considering the following auxiliary problem inspired from \cite{Chitour2019}
\begin{equation}\label{probd}
\left\lbrace
\begin{array}{cccc}
y_{tt} -y_{xx} + a(x)\theta(t,x)y_t=0 & \hbox{for } (t,x) \in  \mathbb{R_+} \times (0,1), \\
y(t,0)= y(t,1)=0   &  t\geq 0 ,\\
y(0,\cdot)= y_0\ ,\ y_t(0,\cdot)=y_1,
\end{array}
\right.
\end{equation}
where $a$ satisfies Hypothesis $\mathbf{ (H_1)}$ and $\theta$ satisfies:\\ \\
$\mathbf{ (H_3)}$ $\theta:\ \mathbb{R_+}\times [0,1] \rightarrow \mathbb{R}$ is a non-negative continuous function such that 
\begin{align}
\exists\ \theta_1,\ \theta_2 >0, \ \theta_1\leq \theta(t,x)\leq  \theta_2 \ \ \forall \ (t,x)\in \mathbb{R_+}\times [0,1].
\end{align}
\begin{Remark}
The well-posedness of Problem \eqref{probd} can be treated the same way as the linear problem was treated in \cite{kafnemer2022}.
\end{Remark}
We define the Riemann invariants for all $(t,x) \in \mathbb{R_+}\times (0,1)$ by
\begin{align}\label{riemann}
\bar{\rho}(t,x)= {y_x(t,x)+y_t(t,x)}, \\
\bar{\xi}(t,x)= {y_x(t,x)-y_t(t,x)}.
\end{align}
Along strong solutions of \eqref{probd}, we deduce that
\begin{equation}\label{probwnl}
\left\lbrace
\begin{array}{ll}
\bar{\rho}_t - \bar{\rho}_x  =-\frac{1}{2}a(x)\theta(t,x)(\bar{\rho}-\bar{\xi}) &\text{in }   \mathbb{R_+} \times (0,1), \\
\bar{\xi}_t + \bar{\xi}_x  =\frac{1}{2}a(x)\theta(t,x)(\bar{\rho}-\bar{\xi}) &\text{in }   \mathbb{R_+} \times (0,1), \\
\bar{\rho}(t,0)-\bar{\xi}(t,0)= \bar{\rho}(t,1)-\bar{\xi}(t,1)=0   &  \forall t \in \mathbb{R_+} ,\\
\bar{\rho}_0:=\bar{\rho}(0,.)= {y_0'+ y_1}\ ,\ \bar{\xi}_0 :=\bar{\xi}(0,.)={y_0'-y_1},
\end{array}
\right.
\end{equation}
with $\left(\bar{\rho}_0, \bar{\xi}_0 \right) \in W^{1,p}(0,1)\times W^{1,p}(0,1)$.
\\ \\
We consider the $p$th-energy $E_p(y)$ of a solution $y$, defined on
$\mathbb{R_+}$  by
\begin{align}
E_p(y)(t)= \frac{1}{p} \int_0^1 (|{\bar{\rho}}|^p+|{\bar{\xi}}|^p) dx.
\end{align}
We deduce from Proposition \ref{prop1} the following proposition.
\begin{Proposition}\label{prop1nl}
	Let $p \in [1, \infty)$ and suppose that a strong solution $y$ of \eqref{probd} exists and is defined on
	a non trivial interval $I \subset \mathbb{R_+}$ containing $0$, for some initial conditions $(y_0,y_1)\in
	Y_p$. For $t\in I$, define
	\begin{align}\label{phidefnl}
	\Phi (t) := \int_0^1 [{\cal{F}}(\bar{\rho})+{\cal{F}}(\bar{\xi})] dx,
	\end{align}
	where $\bar{\rho}$ and $\bar{ \xi}$ are defined in \eqref{riemann} and ${\cal{F}} $ is a $C^1$ convex function. Then $\Phi$ is well defined
	for $t\in I$ and satisfies
	\begin{align}\label{phidecnl}
	\frac{d}{dt} \Phi(t) =-\frac{1}{2}\int_0^1 a(x)\theta(t,x) (\bar{\rho} -\bar{ \xi})({\cal{F}}'(\bar{ \rho})- {\cal{F}}'(\bar{ \xi})) dx \leq 0.
	\end{align}
\end{Proposition}
\textbf{\emph{Proof.}} 	
The proof is the same as the proof of Proposition \ref{prop1} where we replace $g(z_t)$ by $\theta z_t$.
\begin{small}
	\begin{flushright}
		$\blacksquare$
	\end{flushright}
\end{small}
Similarly to Corollary \eqref{cor:decrease}, we have the following.
	\begin{Corollary}\label{cor:decreasenl}
If $(y_0,y_1) \in Y_p$, then we have along strong solutions that for $t\geq 0$,
		\begin{align}\label{E_p'nl}
		E_p(y)'(t) = -\frac{1}{2}\int_0^1 a(x)\theta(t,x)(\bar{ \rho}-\bar{ \xi})\left(\lfloor \bar{ \rho} \rceil^{p-1}-
		\lfloor\bar{ \xi} \rceil^{p-1} \right) dx,
		\end{align}
			where $\bar{\rho}$ and $\bar{ \xi}$ are defined in 
			\eqref{riemann}.
	Moreover,	 for $(y_0,y_1) \in X_p$, suppose that a weak solution $y$ of \eqref{probd}
		exists on $\mathbb{R_+}$. Then the energy $t\longmapsto E_p(y)(t)$ is non-increasing. 
	\end{Corollary}
	The main result of this section is given below.
\begin{Theorem}\label{propaux}
	Fix $ p \in ]1, \infty)$ and suppose that Hypotheses $\mathbf{ (H_1)}$ and $\mathbf{ (H_3)}$ are satisfied. Then for every $(y_0,y_1) \in X_p$, the solution of \eqref{probd} is exponentially stable. 
\end{Theorem}
\ 
\\ 
To prove Theorem \ref{propaux}, we are going to follow the same steps of the proof of \cite[Theorem 4.3]{kafnemer2022}, all the computations remains the same with the only difference that $a(x)$ is now replaced by $a(x)\theta(t,x)$. We have used three multipliers to treat the linear case in \cite{kafnemer2022}. The multipliers were slightly different in the case where $1<p<2$ than the case where in the case $p\geq 2$ but in both cases we have the same potential occurrences of  $a(x)\theta(t,x)$. As a result, we will only give the sketch of the proof in the case where $p\geq2$ and the case where $1<p<2$ will be treated similarly. The sketch of proof that we provide lists the parts of the proof where $\theta(t,x)$ occurs and how they are easily handled.
\\ 
\\
It is important to note that just like in \cite[Theorem 4.3]{kafnemer2022}, it is enough to prove Theorem \ref{propaux} for strong solutions and then extend the result for weak solutions by a density argument.
\subsection{Case $\mathbf{p\geq2}$}
 The theorem for strong solutions and for $p\geq 2$ follows directly from the next proposition by using Gronwall Lemma.
\begin{Proposition}\label{propgeq2}
	Fix $2\leq p<+\infty$ and suppose that the hypotheses of Theorem \ref{propaux} are satisfied. Then there exist positive
	constants $C$ and $C_p$ such that, for every $(z_0, z_1)\in Y_p,$ it holds the following energy estimate:
	\begin{align}\label{energest1}
\forall 0\leq S\leq T, \quad \int_S^T E_p(t)\, dt \leq C C_p E_p(S).
	\end{align}
\end{Proposition}
To prove this key proposition, we divide the proof into steps, the result of each step is given by a key lemma and is obtained by using a specific multiplier. Before we announce the lemmas, we introduce some useful functions.
\\ \\
Let the function $f$ defined by
\begin{equation}\label{fdef}
f(s)=\lfloor s \rceil^{p-1},\qquad \forall\ s\in\mathbb{R}.
\end{equation}
and the function $F(s)= \int_0^s f(\tau) d\tau$, we have that
\begin{equation}\label{Fdef}
F(s) = \frac{|s|^p}{p},\quad F'=f,\qquad  f'(s)=(p-1)\vert s\vert^{p-2}.
\end{equation}
The multipliers that were used in \cite[Section 4.1]{kafnemer2022} in the case $p\geq 2$ are the following: 
\begin{itemize}
	\item[$(m1)$] $x\mapsto x\psi(x) f(\bar{ \rho}(t,x))$ and $x\mapsto x\psi(x)f(\bar{ \xi}(t,x))$ for every $t\geq 0$;
	\item[$(m2)$] $x\mapsto \phi(x) f'(\bar{ \rho}(t,x)) y(t,x)$ and $x\mapsto \phi(x) f'(\bar{ \xi}(t,x)) y(t,x)$
	for every $t\geq 0$;
\end{itemize}
\begin{itemize}
	\item[$(m3)$] $x\mapsto v(t,x)$ for every $t\geq 0$, where $v$ is the solution of the following elliptic problem defined for every $t\geq 0$:
	\begin{equation}\label{elliptnl}
	\left\lbrace
	\begin{array}{ll}
	v_{xx} = \beta f(y)  & x\in (0,1) , \\
	v(0)=v(1)=0,&
	\end{array}
	\right.
	\end{equation}
\end{itemize}
where $\psi, \ \phi, \ \beta$ are the localization functions defined as
\begin{equation}\label{psidef}
\begin{cases}
0\leq \psi \leq 1,
\\
\psi = 0\ \textrm{on}\ Q_0 ,
\\
\psi=1\ \textrm{on}\ (0,1)\setminus Q_1,
\end{cases}
\quad
\begin{cases}
0\leq \phi \leq 1,
\\
\phi = 1\ \hbox{on}\ Q_1 ,
\\
\phi=0\ \hbox{on}\ (0,1)\setminus Q_2,
\end{cases}
\quad \begin{cases}
0 \leq \beta \leq 1,
\\
\beta = 1\ \ \hbox{on}\ \ Q_2 \cap (0,1) ,
\\
\beta=0\ \ \hbox{on}\  \mathbb{R} \setminus \omega.
\end{cases}
\end{equation}
The tests functions are represented in the following figure that was taken from \cite{kafnemer2022},
\begin{center}
	\includegraphics[width=14cm]{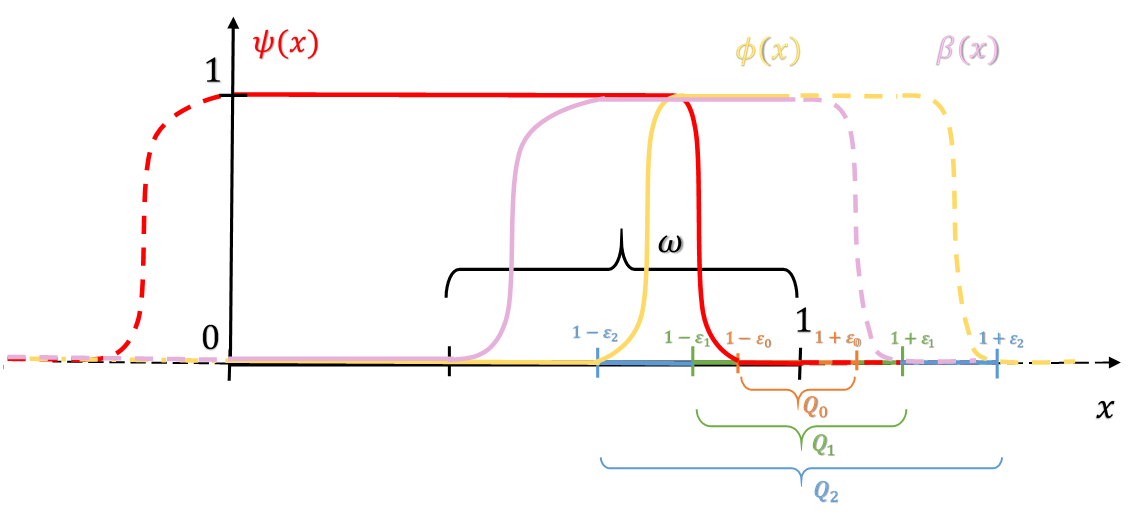}
\end{center}
where the intervals $Q_i$ are defined by small enough ordered positive constants $\epsilon_i$, $0\leq i\leq 2$, accoring to the picture.
\begin{Lemma}\textbf{(First set of multipliers)} \\
Under the hypothesis of Proposition \ref{propaux}, we have for all $0\leq S\leq T$ the following estimate: 
\begin{align}\label{Tnl}
\int_S^T E_p(y)(t) dt \leq CC_p\,E_p(y)(S) +
C\underbrace{\int_S^T \int_{Q_1\cap (0,1)} (F(\bar{ \rho})+F(\bar{ \xi})) \, dx\, dt}_{\mathbf{S_4}}.
\end{align}
\end{Lemma}
\textbf{\emph{Proof:}}
Multiplying the first equation of \eqref{probwnl} by $x\, \psi\, f(\bar{ \rho})$ and then by $x\, \psi\, f(\bar{ \xi})$   and integrating over
$[S,T]\times [0,1]$, we obtain after doing the same manipulations that led to \cite[Equation (4.14)]{kafnemer2022}, but instead we obtain 
\begin{align}\label{T123nl}
\int_S^T E_p(y)(t) dt \leq& \underbrace{ \int_S^T \int_{Q_1\cap (0,1)} |\left(1-(x\, \psi)_x \right)| \left(F(\bar{ \rho})+F(\bar{ \xi})\right) \, dx\, dt}_{\mathbf{S_1}} \notag
\\
&+ \underbrace{\int_0^1 |x\,\psi | \,\left| \left[F(\bar{ \rho})- F(\bar{ \xi})\right]_S^T \right|dx}_{\mathbf{S_2}} \notag \\&
\underbrace{+\frac{1}{2}\int_S^T \int_0^1 |a(x)\theta(t,x)x\psi|\left|(f(\bar{ \rho})+f(\bar{ \xi}))\right||\bar{ \rho}-\bar{ \xi}|\, dx\, dt}_{\mathbf{S_3}} .
\end{align}
The quantities $\mathbf{S_1}$ and $\mathbf{S_2}$ are denoted the same way as $\mathbf{S_1}$ and $\mathbf{S_2}$ from the proof of \cite[Lemma 4.6]{kafnemer2022} and treated the same way to obtain \cite[Equation (4.15)]{kafnemer2022} and \cite[Equation (4.16)]{kafnemer2022}. As for the estimation of $\mathbf{S_3}$, also denoted as $\mathbf{S_3}$ from the proof of \cite[Lemma 4.6]{kafnemer2022}, nothing radically changes, we just follow the same computations as the one that led to \cite[Equation (4.22)]{kafnemer2022} with taking in consideration $\theta(t,x)$,
which gives the same estimation despite the presence of $\theta(t,x)$.
This allows us to combine the estimations of $\mathbf{S_1}$, $\mathbf{S_2}$ and $\mathbf{S_3}$ to obtain the same main result \cite[Equation (4.6)]{kafnemer2022} of the first set of multipliers, which is given in our case by \eqref{Tnl}.
\begin{flushright}
	\begin{small}
		$\blacksquare$
	\end{small}
\end{flushright}
To estimate $\mathbf{S_4}$, we use the second set of multipliers.
\begin{Lemma}\textbf{(Second set of multipliers)} \\
Under the hypotheses of Proposition \ref{propaux}, we have for all $0\leq S\leq T$ the following estimate:
{\begin{align}\label{T3nl}
	\mathbf{S_{4}} \leq  C\frac{C_p}{\eta^p}\underbrace{\int_S^T\int_{Q_2 \cap
			(0,1)} |y|^p \, dx\, dt}_{\mathbf{T_5}}+ C C_p\eta^q
	\int_S^T E_p(y)(t)\,dt +CC_pE_p(y)(S),
	\end{align}}
where $\eta>0$ is arbitrary and $q=\frac{p}{p-1}$.
\end{Lemma}
\textbf{\emph{Proof:}}
Multiplying the first equation of \eqref{probwnl} by $\phi f'(\bar{ \rho})y$ and then by $\phi f'(\bar{ \xi})y$ and integrating over
$[S,T]\times [0,1]$, we obtain after doing the same manipulations that led to \cite[Equation (4.35)]{kafnemer2022}, but instead we obtain 
\begin{align}\label{mul2'nl}
\mathbf{S_{4}}  &\leq C\underbrace{\int_S^T\int_{Q_2 \cap (0,1)} |y| \left(|f(\bar{ \rho})| + |f(\bar{ \xi})| \right) \, dx\, dt}
_{\mathbf{T_1}}+ C_p\underbrace{\left|\left[ \int_0^1\left(f(\bar{ \rho})-f(\bar{ \xi})\right) y\, dx\right]_S^T\right| }
_{\mathbf{T_2}} \notag \\
&\ \ \ \ \ \ \ +C_p\underbrace{\int_S^T\int_{Q_2 \cap (0,1)} |\left(f'(\bar{ \rho}) + f'(\bar{ \xi}) \right) y\, a(x)\theta(t,x)(\bar{ \rho} - \bar{ \xi})| \, dx\, dt}_{\mathbf{T_3}}  
\notag \\
&\ \ \ \ \ \ \ \ \ \ \ \ \ \ \ \ \ +C_p \underbrace{\int_S^T\int_0^1 \left|\phi (\bar{ \rho}-\bar{ \xi}) \left(f(\bar{ \rho}) - f(\bar{ \xi})\right)\right| \, dx\, dt}_{\mathbf{T_4}}.
\end{align}
${\mathbf{T_1}}$ and ${\mathbf{T_2}}$ are denoted the same way as $\mathbf{T_1}$ and $\mathbf{T_2}$ in the proof of \cite[Lemma 4.7]{kafnemer2022} and treated the same way to obtain \cite[Equation (4.37)]{kafnemer2022} and \cite[Equation (4.40)]{kafnemer2022}. As for estimating ${\mathbf{T_3}}$, also denoted by $\mathbf{T_3}$ in the proof of \cite[Lemma 4.7]{kafnemer2022}, we just use the fact that $ \theta_1\leq \theta(t,x)\leq \theta_2$ to obtain
\begin{equation}
\vert a(x) \theta(t,x)\left(f'(\bar{ \rho}) + f'(\bar{ \xi}) \right) (\bar{ \rho} - \bar{ \xi})\vert\leq C_p\left(\vert f(\bar{ \rho})\vert+\vert f(\bar{ \xi})\vert\right),
\end{equation} 
which implies that ${\mathbf{T_3}}\leq {\mathbf{T_1}} $ and hence it has the same estimation as $\mathbf{T_1}$. 
\\ \\
Left to estimate ${\mathbf{T_4}}$ which is also denoted by $\mathbf{T_4}$ in the proof of \cite[Lemma 4.7]{kafnemer2022}. We follow the same steps that led to \cite[Equation (4.42)]{kafnemer2022} but we use the fact that $|\phi(x)| \leq Ca(x)\theta(t,x)$ for all $t\in \mathbb{R_+}$ and $x\in [0,1]$. We obtain the same estimation as \cite[Equation (4.42)]{kafnemer2022}.\\ \\
We combine the estimations of $\mathbf{T_1}$, $\mathbf{T_2}$ and $\mathbf{T_3}$, we obtain the main result of the second set of multipliers, which is the same estimation \cite[Equation (4.25)]{kafnemer2022} given in our case by
\eqref{T3nl}. 
\begin{flushright}
	\begin{small}
		$\blacksquare$
	\end{small}
\end{flushright}
Left to estimate now ${\mathbf{T_5}}$. To do so, we use one last multiplier.
\begin{Lemma}\textbf{(Third multiplier)} \\
	Under the hypotheses of Proposition \ref{propaux}, we have for all $0\leq S\leq T$ the following estimate:
	\begin{align} \label{S1nl}
	\underbrace{\int_S^T \int_{Q_2\cap(0,1)}|y|^p \, \, dx\, dt}_{\mathbf{T_5}} \leq  CC_p\left(\eta
	\int_S^T E_p(y)(t) \, dt+\frac1{\eta^r} E_p(y)(S)\right).
	\end{align}
\end{Lemma}
\textbf{\emph{Proof:}}
First, one should note that
\cite[Lemma 4.8]{kafnemer2022} remains valid here and it gives an estimation of the $L^q$ norms of $v$ and $v_t$, where $v$ is defined in \eqref{elliptnl}.\\ \\
Multiplying the first equation of \eqref{probwnl} by $v$ and integrating over
$[S,T]\times [0,1]$, we obtain after doing the same manipulations that led to \cite[Equation (4.58)]{kafnemer2022}, bu instead we obtain 
\begin{align} \label{thirdmulnl}
2{\mathbf{T_5}} &\leq  \underbrace{ \left|\left[ \int_0^1 v (\bar{ \rho}-\bar{ \xi}) dx \right]_S^T\right|}_{\mathbf{V_1}} + \underbrace{ \int_S^T \int_0^1 |v_t| |(\bar{ \xi} - \bar{ \rho})|\, \, dx\, dt}_{\mathbf{V_2}} \notag\\
& +\underbrace{ \int_S^T \int_0^1 |v a(x)\theta(t,x)(\bar{ \rho}-\bar{\xi})|\, \, dx\, dt}_{\mathbf{V_3}}.
\end{align}
The quantities ${\mathbf{V_1}}$, ${\mathbf{V_2}}$ and ${\mathbf{V_3}}$ are denoted the same way as ${\mathbf{V_1}}$, ${\mathbf{V_2}}$ and ${\mathbf{V_3}}$ from the proof of \cite[Lemma 4.9]{kafnemer2022}.
The quantity ${\mathbf{V_1}}$ is treated the same way as\cite[Equation (4.60)]{kafnemer2022}. As for ${\mathbf{V_3}}$, we follow the same steps that led to \cite[Equation (4.64)]{kafnemer2022} and we use the fact that  $\beta\leq C\, a(x) \theta(t,x)$ in \cite[Equation (4.62)]{kafnemer2022}. Finally for ${\mathbf{V_3}}$, nothing changes and we obtain \cite[Equation (4.68)]{kafnemer2022} despite the presence of $\theta(t,x)$ in this term. We combine the three estimations to obtain the main result of the third multiplier, which is the same estimation of ${\mathbf{T_5}}$ given by \cite[Equation (4.49)]{kafnemer2022}, given in our case by \eqref{S1nl}.
\begin{flushright}
	\begin{small}
$\blacksquare$
\end{small}
\end{flushright}
\textbf{\emph{Proof of Proposition \ref{propgeq2}:}} \\ \\
Finally, as in \cite[Section 4.1.4]{kafnemer2022}, by combining the results of the three multipliers \eqref{Tnl}, \eqref{T3nl} and \eqref{S1nl} and choosing $\eta$ properly, we obtain the energy estimate given by \eqref{energest1}.
\begin{flushright}
	\begin{small}
		$\blacksquare$
	\end{small}
\end{flushright}
\subsection{Case $\mathbf{1<p<2}$}
We remind the reader that due to the presence of the power $p-2$ in the second set of  multipliers in the case $p\geq 2$, it is not possible to use them directly in the case $1<p<2$. Therefore, we have to modify the functions $f$ and $F$ just like it was done in \cite[Section 4.2]{kafnemer2022}. We
consider, for $p\in (1,2)$, the functions $g$ and $G$ defined on $\mathbb{R}$, by
\begin{align}
g(y)&= (p-1)\int_0^y (|s| + 1)^{p-2} \, ds =\textrm{sgn}(y)\left[(|y|+1)^{p-1}-
1\right], \label{gdef}\\
G(y)&=\int_0^y g(s)\,ds =\frac{1}{p}\left[(|y|+1)^p-1\right]-|y|.
\label{Gdef}
\end{align}
We also modify
the energy $E_p$ by considering,  for every
$t\in \mathbb{R_+}$ and every solution of \eqref{probnl}, the function ${\cal{E}}_p$ defined by
\begin{align}\label{Eepsdef}
{\cal{E}}_p(t)= \int_0^1 \left( G(\rho) +G(\xi) \right) \, dx.
\end{align}
The proof of Theorem \ref{propaux} is a result of the following proposition by using the argument given by \cite[Section 4.2.4]{kafnemer2022}.
\begin{Proposition}\label{prop<2}
Fix $p\in (1,2)$ and suppose that the hypotheses of Theorem \ref{propaux} are satisfied. Then there exist positive
	constants $C$ and $C_p$ such that, for every $(z_0,z_1)\in Y_p$ verifying
\begin{align}
E_p(0)\leq 1,
\end{align}
	we have the following energy estimate:
	\begin{align}
	\forall 0\leq S\leq T, \ \ \int_S^T {\cal{E}}_p \, dt\leq  CC_p{\cal{E}}_p,
	\end{align}
	where ${\cal{E}}_p$ is defined in \eqref{Eepsdef}.
\end{Proposition}
\textbf{\emph{Proof:}}
Just like the case $p\geq 2$, we follow the exact same steps in this case as the case ${1<p<2}$ from \cite[Section 4.2]{kafnemer2022}, by using the multipliers 
\begin{description}
	\item[$(\bar{m}1)$] $x\mapsto x\psi(x) g(\bar{ \rho}(t,x))$ and $x\mapsto x\psi(x)g(\bar{ \xi}(t,x))$ for every $t\geq 0$;
	\item[$(\bar{m}2)$] $x\mapsto \phi(x) g'(\bar{ \rho}(t,x)) y(t,x)$ and $x\mapsto \phi(x) g'(\bar{ \xi}(t,x)) y(t,x)$
	for every $t\geq 0$;
	\item[$(\bar{m}3)$] $x\mapsto v(t,x)$ for every $t\geq 0$, where $v$ is the solution of the following elliptic problem defined for every $t\geq 0$:
	\begin{equation*}\label{elliptnl<2}
	\left\lbrace
	\begin{array}{ll}
	v_{xx} = \beta g(y)  & x\in (0,1) , \\
	v(0)=v(1)=0,&
	\end{array}
	\right.
	\end{equation*}
\end{description}
where $\psi, \ \phi, \ \beta$ are the localization functions defined in \eqref{psidef} and $g$ is the function defined in \eqref{gdef}.\\ \\
 We take care of the presence of $\theta$ just like we did in the previous proof since we have technically the same occurrences of  $\theta$ in the cases $p\geq 2$ and $1<p<2$. We obtain the exponential decay similarly.
\begin{flushright}
	\begin{small}
		$\blacksquare$
	\end{small}
\end{flushright}
Now that we know that Problem \eqref{probd} is exponentially stable for all $1<p<+\infty$, we can now consider Problem \eqref{probnl}.
\subsection{Asymptotic stability of the nonlinear problem}
We conclude in this section by using what precedes, the stability of the nonlinear problem. However, to be able to finally conclude we need to state and prove some key lemmas first.
\begin{Lemma}
	Define the function the continuous function $\nu$ defined for all $x\in \mathbb{R}$ as 
	\begin{equation*}
	\left\lbrace
	\begin{array}{ccc}
	\nu(x)=\frac{g(x)}{x}& \ \ \ \hbox{ for } x \in  \mathbb{R^*}, \\
	\nu(0)=g'(0),& 
	\end{array}
	\right.
	\end{equation*}
	where $g$ is defined in Hypothesis $\mathbf{(H_2)}$, then $\nu$ has the following properties: 
	\begin{itemize}
		\item $\nu(x)>0$ for all $x\neq 0$.
		\item $\forall M>0$, $\exists \nu_1(M), \nu_2(M)>0$ such that $\forall |x|\leq M$,\ \ $ \nu_1(M)\leq \nu(x) \leq \nu_2(M)$. 
	\end{itemize}
\end{Lemma}
\textbf{\emph{Proof:}} The proof is standard and is a direct result of Hypothesis $\mathbf{(H_2)}$. The first item is a direct result of \eqref{gxx}. The second item is a result of $g$ being $C^1$, $g(0)=0$ and $g'(0) \neq 0$.
\begin{small}
	\begin{flushright}
		$	\blacksquare$
	\end{flushright}
\end{small}
We next prove the following regularity lemma which is a generalization of \cite[Lemma 2]{Martinez2000}.
\begin{Lemma} \label{boundu'} 
	Assume $\mathbf{(H_1)}$ and $\mathbf{(H_2)}$ satisfied and $(z_0,z_1) \in Y_p$ with $1<p<+\infty$, then the following estimate holds true for all $t\geq 0$ and for a constant $C_p(z_0,z_1)>0$ that depends on the norm of initial conditions in $Y_p$
	\begin{align}
	||z_t(t,\cdot)||_{W^{1,p}(0,1)} \leq C_p(z_0,z_1). \label{boundrx}
	\end{align}
\end{Lemma}
\emph{\textbf{Proof:}}
The norm estimate \eqref{boundrx} is a direct result of \eqref{w}.
The lemma is then concluded with
\begin{align}
C_p(z_0,z_1)=(CpE_p(w) (0))^{\frac{1}{p}},
\end{align}
where $w=z_t$.
\begin{small}
	\begin{flushright}
		$	\blacksquare$
	\end{flushright}
\end{small}
\textbf{\emph{Proof of Theorem~\ref{stability}}}: \\ \\
Consider an arbitrary strong solution  $z$ of Problem \eqref{probnl}. Using the result of Lemma~\ref{boundu'} which is $z_t \in W^{1,p}(0,1)$ with $||z_t(t,\cdot)||_{W^{1,p}(0,1)}^p \leq C_p(z_0,z_1)$ for all $t\in \mathbb{R_+}$ we deduce using the continuous embedding $W^{1,p}(0,1) \subset L^\infty (0,1),$ with a constant that does not depend on $t$  (see \cite[Theorem VIII.7]{Brezis}) that
\begin{align}\label{ztinf}
||z_t(t,\cdot)||_{L^\infty(0,1)} \leq C ||z_t(t,\cdot)||_{W^{1,p}(0,1)} \leq C_p(z_0,z_1) \ \ \forall \ t\in \mathbb{R_+}.
\end{align}
Since $g$ satisfies $\mathbf{ (H_2)}$ and $z_t$ satisfies \eqref{ztinf}, we deduce the existence of two constants $$C^1_p(z_0,z_1), C^2_p(z_0,z_1)>0,$$ that depends on $p$ and on the $Y_p$-norm of the initial data only such that for all $t\geq 0$, and $x\in(0,1)$, it holds 
\begin{align}\label{g(z_t)}
C^1_p(z_0,z_1) \leq \nu(z_t) \leq C^2_p(z_0,z_1).
\end{align}
We consider the time-varying linear problem given by
\begin{equation}\label{probww}
\left\lbrace
\begin{array}{cccc}
y_{tt} -y_{xx} + a(x)\nu(z_t)y_t=0 & \hbox{for } (t,x) \in  \mathbb{R_+} \times (0,1), \\
y(t,0)= y(t,1)=0   &  t\geq 0 ,\\
y(0,\cdot)= y_0\ ,\ y_t(0,\cdot)=y_1,
\end{array}
\right.
\end{equation}
which is nothing else but the auxiliary problem \eqref{probd} with $\theta(t,x)=\nu(z_t)$ since $\nu(z_t)$ is seen as a function of $t$ and $x$. Moreover, $\nu(z_t)$ satisfies Hypothesis $\mathbf{ (H_3)}$ with $\theta_1=C^1_p(z_0,z_1)$ and $\theta_2=C^2_p(z_0,z_1)$. 
\\
\\
By Theorem \ref{propaux}, the energy $E_p(y)$ decays exponentially to zero along the solutions of \eqref{probww}. The key and trivial remark is that the strong solution $z$ of Problem \eqref{probnl} considered at the beginning of the argument 
 is the solution of \eqref{probww} with initial data $(z_0,z_1)$.
 This concludes the proof of Theorem \ref{stability}.
\begin{small}
	\begin{flushright}
		$	\blacksquare$
	\end{flushright}
\end{small}
 \bibliographystyle{plain}
\bibliography{biblio}
\end{document}